\documentclass{amsart}
\makeatletter
\@namedef{subjclassname@2020}{%
  \textup{2020} Mathematics Subject Classification}
\makeatother
\usepackage{latexsym}
\usepackage{amssymb}
\usepackage{graphicx}
\usepackage{wrapfig}
\usepackage{amsthm}
\usepackage{float}
\usepackage{epsfig}
\usepackage{multirow}
\usepackage{xcolor}
\usepackage{caption}
\usepackage[toc,page]{appendix}
\usepackage{old-arrows} 
\usepackage{caption}
\usepackage{subcaption}
\usepackage{lmodern, mathtools}
\usepackage{mathtools}
\usepackage{hyperref}
\RequirePackage{color}
\usepackage{tikz}
\usepackage{amscd,enumerate}
\usepackage{amsfonts}
\hypersetup{ colorlinks,
linkcolor=blue,
filecolor=green,
urlcolor=blue,
citecolor=blue }
\usepackage{helvet}
\usepackage{bigstrut}
\usepackage{float, caption, booktabs, cellspace}
\usepackage{calrsfs}
\DeclareMathAlphabet{\pazocal}{OMS}{zplm}{m}{n}

\makeatletter
\@addtoreset{subfigure}{framenumber}
\makeatother

\usepackage{etoolbox}
\AtBeginEnvironment{figure}{\setcounter{subfigure}{0}}

\newtheorem*{theo}{Theorem}
\newtheorem*{conj}{Conjecture}
\newtheorem*{prop}{Proposition}
\newtheorem*{coro}{Corollary}
\newtheorem*{lem}{Lemma}

\input xygraph
\input xy
\xyoption{all}

\usepackage{array}
\newcolumntype{L}[1]{>{\raggedright\let\newline\\\arraybackslash\hspace{0pt}}m{#1}}
\newcolumntype{C}[1]{>{\centering\let\newline\\\arraybackslash\hspace{0pt}}m{#1}}
\newcolumntype{R}[1]{>{\raggedleft\let\newline\\\arraybackslash\hspace{0pt}}m{#1}}

\usepackage{multicol}

\def\g{\gamma}

\def\s{\sigma}

\def\remove#1{}
\newcommand{\p}{{\mathbb{P}}}
\newcommand{\G}{{\mathbb{G}}}

\newcommand{\Z}{{\mathbb{Z}}}

\newcommand{\Q}{{\textit{Q}}}

\title[Homemade Algebraic Geometry]{Homemade Algebraic Geometry.\\ \small  Celebrating Enrique Arrondo's 60th birthday}

\author[S. Marchesi]{Simone Marchesi  $^{1}$ }
\author[A. Tocino]{Alicia Tocino $^{2,\ast}$}

\address{$^{\ast}$ Corresponding author. Email: alicia.tocino@uma.es}

\address{$^{1}$ Facultat de Matemàtiques i Informàtica,
Universitat de Barcelona,
Gran Via de les Corts Catalanes 585,
08007, Barcelona,
Spain}
\address{$^{1}$ Centre de Recerca Matem\`atica Edifici C, Campus Bellaterra, 08193 Bellaterra, Spain}
\email{marchesi@ub.edu}

\address{$^{2}$ Departamento de Matem\'atica Aplicada, E.T.S. Ingenier\'\i a Inform\'atica, Universidad de M\'alaga, Campus de Teatinos s/n. 29071 M\'alaga,   Spain. }
\email{alicia.tocino@uma.es}

\thanks{SM is partially supported by PID2020-113674GB-I00 and by the Spanish State Research Agency, through the Severo Ochoa and Mar\'ia de Maeztu Program for Centers and Units of Excellence in R\&D (CEX2020-001084-M); and he is a member of the GNSAGA group of INDAM. AT is supported by the Spanish Ministerio de Ciencia e Innovaci\'on through project  PID2019-104236GB-I00/AEI/10.13039/501100011033 and  by the Junta de Andaluc\'{\i}a  through the project  FQM-336 with FEDER funds. 
}

\usepackage[foot]{amsaddr}

\begin{document}

\begin{abstract}
In this survey we recognize Enrique Arrondo's contributions over the whole of its career, recalling his professional history and collecting the results of his mathematical production.
\end{abstract}

\maketitle

\section{Introduction}

At July 10th-13th, 2023, we celebrated Enrique Arrondo's 60th birthday. On this occasion we could see and appreciate the professional esteem and the personal affection that so many people, in the mathematical community, have for Enrique.

As a continuation of this celebration, the purpose of this survey is twofold.\\ 
First, we would like to briefly recall, as done in the dedicated talk of the congress, Enrique's academic background and professional history, drawing a global picture of the vast network that he has built over the years and underlying how much he deeply influenced the community (and especially his Ph.D. students) with his unique and \textit{natural} way of doing mathematics. In doing so, we hope we managed to mention most of the people that played a role in his career and apologize in advance if we have missed someone.\\
Furthermore, we also would like to collect in this work most of Enrique's contributions in algebraic geometry, regarding especially Grassmannian varieties and vector bundles. 

Enrique obtained his Degree in Mathematical Science in the Universidad Complutense de Madrid (with extraordinary bachelor's degree award) in 1985. He obtained his Ph.D. in Mathematical Sciences in 1990 from the Universidad Complutense de Madrid (with extraordinary Ph.D. award), presenting the thesis \textit{Congruencias de rectas en $\mathbb{P}^3$} (Congruences of lines in $\mathbb{P}^3$), with advisor Ignacio Sols. Since 1990, he is a professor at the Universidad Complutense de Madrid.

During this period, he provided an important contribution in the training of young researchers, also by mentoring and teaching in summer schools and training courses. His enthusiasm, mastery of algebraic geometry and, at the same time, eagerness to learn, have been, over the years, an inspiration to many. The following genealogy tree shows all of Enrique's mathematical descendants, most of which are now themselves university professors.

$$\scriptsize{\textbf{Ignacio Sols (1975)} \rightarrow 
    \textbf{Enrique Arrondo (1990)}\begin{cases}
        \textbf{Raquel Díaz (1996) - Universidad Complutense de Madrid}\\
        \textbf{Luca Ugaglia (2001) - Università degli studi di  Palermo}\\
        \textbf{Beatriz Graña (2003) - Universidad de Salamanca}\\
        \textbf{Jose Carlos Sierra (2004) - Universidad Nacional de Educación a Distancia}\\
        \textbf{Jorge Caravantes (2006) - Universidad de Alcalá}\\
        \textbf{Sofía Cobo (2008)}\\
        \textbf{Simone Marchesi (2012) - Universitat de Barcelona} \begin{cases}
            \textbf{Aydee López (2017)}\\
            \textbf{Aline Vilela (2018)}
        \end{cases}\\
        \textbf{Alicia Tocino (2015) - Universidad de Málaga}
    \end{cases}\\}
$$

\bigskip

We now pass to Enrique's contributions in algebraic geometry.\\ 
Being aware of how hard and unfair it is to divide a life's work in sections, we have decided to start by presenting the research developed during his doctoral studies, focused on \textit{congruences} over projective spaces and Grassmannian varieties. Among his collaborators on this topic over the years we would like to mention Marina Bertolini, Sofía Cobo, Beatriz Graña, Mark Gross and Cristina Turrini. Nevertheless, more details on this topic can be found in Section \ref{sec-congruences}.

To continue this division, an efficient way to detect Enrique's mathematical interests is given by considering the topics proposed to his Ph.D. students.

Hence, wee see that another point of interest consists in \textit{vector bundles without intermediate cohomology}. Together with Beatriz Graña, Francesco Malaspina and Alicia Tocino these bundles were investigated on  Grassmannians of lines. Furthermore, this study was extended for Fano $3$-fold and quartic threefolds together with Daniele Faenzi, Laura Costa and Carlo Madonna. Details on this topic can be found in Section \ref{sec-nocohomology}.

Moving forward, the works together with Raffaella Paoletti, Jose Carlos Sierra and Luca Ugaglia provide the next point: \textit{projections of Grassmannians}.  Details on this topic can be found in Section \ref{sec-projections}.

The following item of our list is given by \textit{Steiner and Schwarzenberger bundles}. They were first studied by Enrique for projective spaces and subsequently generalized, first for Grassmannians, with Simone Marchesi, and then for projective varieties, jointly with Simone Marchesi and Helena Soares. Details on this topic can be found in Section \ref{sec-Steiner}.

The next area of interest is to be found in the famous \textit{Hartshorne's conjecture}, the study of which inspired various lines of research. An example of this is the subcanonicity of codimension two submanifolds of $\G(1,4)$, he studied together
with Maria Lucia Fania, as well as the Picard group of low codimension subvarieties, studied with Jorge Caravantes. 
Details on this topic can be found in Section \ref{sec-Hartshorne}.

Finally, Section \ref{sec-various} collects ``miscellaneous'' topics and Section \ref{sec-EnriqueAG} represents an exhibit of Enrique's unique point of view in understanding algebraic geometry.

We cannot finish this introduction without thanking Enrique again, who helped us grow professionally and personally, and whose friendships we hold dear.

\begin{center}
\includegraphics[width=4cm]{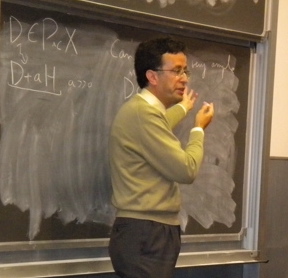}
    \\
    \tiny{Projective Algebraic Geometry in Milano, June, 2009\\
    Courtesy of Antonio Lanteri}
\end{center}

\section{Congruences}\label{sec-congruences}
As we mentioned before, Enrique's Ph.D. thesis is entitled \textit{Congruencias de rectas en $\mathbb{P}^3$}; it has been developed under the supervision of Ignacio Sols and defended in 1990. In their joint published work \cite{Sols1992Congruences},
they mainly study smooth congruences, which are surfaces in the Grassmannian of lines in $\p^3$, that is, $\G(1,3)$. They also inquire the parallelism with surfaces in $\p^4$.
They prove that the only indecomposable bundles on $\G(1,3)$, without intermediate cohomology, are line bundles and twists of the spinor bundle (we will briefly explain the concept of vector bundles without intermediate cohomology in the next section). They also describe the Hilbert schemes of all smooth congruences of degree up to nine, enhancing the results obtained previously in \cite{Sols1989} and by Alessandro Verra in \cite{Verra}. The most relevant result is \cite[Theorem 5.1]{Sols1992Congruences}, where they give a complete classification of the smooth congruences that can be obtained as a projection from another surface in $\G(1,4)$ (projections will be another important topic of Enrique's research that will be detailed later on). Although it is an analogue of Severi's theorem for $\p^4$ (see \cite{Severi}) they use a completely different approach. In fact, they obtain five different classes of such smooth congruences geometrically, more concretely by looking at the geometry of lines in $\p^3$. Moreover, with the collaboration of Manuel Pedreira, they prove:
\begin{theo}[\cite{Sols1992Congruences}, Theorem 6.18]
    Except for a finite number of components, each component of the
Hilbert scheme of smooth congruences consists of surfaces of general type.
\end{theo}
\noindent
Indeed, an analogous result for surfaces in $\p^4$ was proven by Geir Ellingsrud and Christian Peskine in \cite{Peskine}.

Let us continue with more congruences, jointly with Mark Gross, \cite{Mark1993}. Contrary to the preferred setting of the classical authors (see \cite{Fano}), they study smooth congruences having a fundamental curve. The fundamental curve of a congruence is formed by its singular points that are points in $\p^3$ with infinitely many lines of the congruence passing through them. 
The number $a$ of lines of a congruence $Y$ passing through a general point of $\p^3$ is called order of $Y$. Dually, the number $b$ of lines of a congruence contained in a general plane of $\p^3$ is called class of $Y$. So, the pair $(a,b)$ is called the bidegree of the congruence.
They provide a more comprehensive classification of smooth congruences in $\G(1,3)$, listing the possible degrees of the fundamental curve and, in case this is not a line, giving all the possible bidegrees of the congruences, that we collect below. 
\begin{theo}[\cite{Mark1993}, Theorem 2.1]
Let $Y$ be a smooth congruence having a curve $C$ in $\p^3$ of singular points. Then one of the following holds:
\begin{enumerate}
    \item[(a)] The curve $C$ is a line.
    \item[(b)] The congruence consists of the bisecants to $C$, which is either a twisted cubic, $Y$ being a $(1,3)$ congruence, or an elliptic quartic and $Y$ is a $(2,6)$ congruence.
    \item[(c)] The congruence $Y$ is a scroll of degree bigger than two and hence is either a $(1,2)$ or a $(2,2)$ (with $C$ being a conic) or a $(3,3)$ congruence and $C$ is a smooth plane cubic. 
    \item[(d)] The curve is a smooth plane cubic and $Y$ is the $(3,6)$ conic bundle over $C$ of Example 1.5 of \cite{Mark1993} (that consists of a concrete surface with an infinite number of singular points of bidegree $(3,6)$).
\end{enumerate}
\end{theo}
\noindent
In addition, they introduce a vector bundle construction for various smooth congruences in Section 3. In particular, in Example 3.7, they give a smooth $(5,8)$  congruence, thanks to an observation of Peskine, being the first example found so far. 

A generalization of these arguments comes about through his collaboration with Marina Bertolini and Cristina Turrini.\\
As a first step,  the three of them study, in \cite{BertoliniTurrini1994Classif}, the congruences of lines, which are defined as $(n-1)$-dimensional subvarieties of $\G(1,n)$. Again, the fundamental curve of the congruence is a curve $C\subset \p^n$ which meets all lines of the congruence. They give a classification of all smooth congruences having a fundamental curve $C$, obtaining two cases. Firstly, when $C$ is a line (Theorem 1), there are infinitely many families of these congruences. Secondly, if ${\rm deg}(C)\geq 2$ (Theorem 2), there are finitely many of such families.
They continue these investigations in \cite{BertoliniTurrini1998SmallDegree}, studying congruences of small degree in $\mathbb{G}(1,4)$. A classification of all smooth threefolds in $\G(1,4)$ is provided in terms of the bidegree $(a,b)$ with $a=0$ or $b\leq 2$ and collected in Lemma 5.1, Lemma 5.2, Lemma 5.3 and Lemma 5.4. Furthermore, all the possible numerical invariants of smooth threefolds in $\G(1,4)$ of degree less than or equal to $10$ are also given. In Table 1 they collect the list of congruences of degree $d\leq 8$, developing all the details in Sections 8, 9 and 10. In Table 2 and Table 3 they provide a maximal list of congruences of degree $9$ and $10$
(which is far from being effective due the presence of unknown cases).  In order to obtain these results, they use different classification results of varieties of small degree, in particular the one given by Paltin Ionescu  in \cite{Ionescu} and the one given by Maria Lucia Fania and Elvira Laura Livorni in \cite{Livorni}. \\
They also continue studying quadric bundle congruences in $\G(1,n)$ with $n\geq 4$ in \cite{BertoliniTurrini2000Grass}. A quadric bundle congruences in $\G(1,n)$ is a quadric fibrations embedded in $\G(1,n)$ with dimension equal to codimension.
They describe all the possible smooth congruences in $\G(1,n)$, for $n\geq 4$, which have a quadric bundle structure over a curve. In particular, Section 3 is devoted to the case $n=6$, Section 4 to the case $n=5$, Section 6 details the case $n=4$ and Section 8 collects all the information in a table. The main tool they use is Castelnuovo's bound for the genus of projective curves as well as a generalization for curves in an arbitrary Grassmannian variety (which was obtained by Luis Giraldo in \cite{Giraldo}).  Much more articles have been published with Marina Bertolini and Cristina Turrini concerning congruences (see \cite{BertoliniTurrini2001Surfaces}, 
\cite{BertoliniTurrini2005Loci} and \cite{BertoliniTurrini2011}). 
Among them, we would like to highlight the following result.
\begin{theo}[\cite{BertoliniTurrini2005Loci}, Theorem 5.1] 
The only smooth congruences $X$ of the trisecant lines to a surface $S$ in $\p^3$ (with at most ordinary singularities) are those listed in the following:
\begin{enumerate}
    \item[(i)] the congruence of trisecants to $S_1$, which is the hyperplane section of $\G(1,4)$ and has bidegree $(1,2)$ and sectional genus $1$. ($S_1$ is the projected Veronese surface)
    \item[(ii)] the congruence of trisecants to $S_3$, which is the dependency locus of two sections of $Q^2$ and has bidegree $(0,2)$ and sectional genus $0$. ($S_3=\text{Bl}_{q,p_1,\ldots,p_7}(\p^2)$ is the Castelnuovo surface of degree $5$, that is, the blowingup of $\p^2$ in eight points embedded by plane quartics with a given double and other seven base points and $Q$ is the quotient bundle of rank $2$ of $\G(1,4)$)
    \item[(iii)] the congruence of trisecants to $S_4$, which is the zero locus of section of $Sym^2Q$ and has bidegree $(0,4)$ and sectional genus $1$. ($S_4=V(2,3)$ is the smooth complete intersection of a quadric and a cubic hypersurface)
    \item[(iv)] the congruence of trisecants to $S_5$, which is the dependency locus of four sections of $Q^3$ and has bidegree $(1,8)$ and sectional genus $10$. ($S_5=\text{Bl}_{p_1,\ldots,p_{10}}(\p^2)$ is the Bordiga surface of degree $6$, embedded by plane quartics with ten base points)
    \item[(v)] the congruence of trisecants to $S_6$, which is the dependency locus of three sections of $Q\oplus Sym^2Q$ and has bidegree $(2,15)$ and sectional genus $33$. ($S_6=\text{Bl}_p(K3)$ is the inner projection of the complete intersection $V(2,2,2)$ of three quadric hypersurfaces in $\p^5$ from a point in it)
    \item[(vi)] the congruence of trisecants to $S_7$, which is the dependency locus of two sections of $Sym^3Q$ and has bidegree $(6,42)$ and sectional genus $181$. ($S_7=V(3,3)$ is the smooth complete intersection of two cubic hypersurfaces)
\end{enumerate}
\end{theo}

Furthermore, Enrique studies on his own line congruences of low order, see \cite{Arrondo2002LineCongr}. As for $\p^3$, the order is defined as the number of lines of the family passing through a general point of $\p^n$. More concretely, Proposition 2.1 characterizes line congruences of order 0 in $\p^4$, Theorem 2.1 gives a complete list of congruences of order one in $\p^3$ and Propositions 3.1, 3.2 and 3.3 provide some results concerning congruences of order two in $\p^3$.

Let us keep on going with more congruences, in this case provided together with Beatriz Graña and Sofía Cobo, both Ph.D. students of Enrique. Beatriz defended her thesis in 2003, entitled \textit{Escisión de fibrados en $\G(1,4)$ y sus variedades}, in which congruences on $\G(1,4)$ with split universal quotient bundle are studied. 
Moreover, Enrique and Beatriz give in \cite[Theorem 2.8]{BeatrizGrana2006} a complete classification of the smooth threefolds in $\G(1,4)$, in the case that the restriction of the universal quotient bundle $Q$ is a direct sum of two line bundles. Sofía Cobo obtained her Ph.D. in 2008, presenting the work \textit{Estabilidad del fibrado universal restringido a congruencias}. In the corrisponding published paper \cite{Cobo},  Enrique and Sofia try to discover if there is a congruence in $\G(1,3)$ with any preassigned bidegree $(a,b)$ obtaining as main result in this sense Theorem 2.1.

\section{Vector bundles without intermediate cohomology}\label{sec-nocohomology}
Now, it is time for arithmetically Cohen-Macaulay (aCM for short) vector bundles, equivalently, vector bundles without intermediate cohomology.
Our story starts when Geoffrey Horrocks obtains in \cite{Horrocks} a criterion that states that a vector bundle $F$ over $\p^n$ splits as a direct sum of line bundles 
if and only if $H^i(\p^n,F\otimes \mathcal{O}_{\p^n}(t))=0$ for all $t\in\Z$ and $0<i<n$, that is, $F$ does not have intermediate cohomology. In this direction, Horst Knörrer proves that line bundles and twists of the spinor bundles are the only indecomposable aCM
vector bundles over quadrics in \cite{Knorrer}. In a kind of converse result, Buchweitz, Greuel and Schreyer show in \cite{GBS} that the only smooth hypersurfaces in a projective space for which there exists, up to a
twist, a finite number of aCM bundles are the hyperplanes and the quadrics.
Continuing our story, Giorgio Ottaviani generalizes Horrocks' criterion for quadrics and Grassmannians, respectively in \cite{OttavianiGrass2} and \cite{OttavianiGrass}. In the last case, he states that a vector bundle $F$ over $\G(k,n)$ splits if and only if $H^i(\G(k,n),\wedge^{i_1}\Q \otimes\ldots\otimes\wedge^{i_s}\Q \otimes F(t))=0$ for all  $0\leq i_1,\ldots,i_s\leq n-k$, $s\leq k$,  $t\in\Z$ and $0<i<(k+1)(n-k)$ where $\Q$ is the quotient bundle on $\G(k,n)$.
In the same sense as Horrocks' and Ottaviani's criterion, Enrique Arrondo and Beatriz Graña  characterize exactly which are the concrete vector bundles on $\G(1,4)$ without intermediate cohomology in \cite[Theorem 3.3]{BeatrizGrana1999}. Following the Mumford-Castelnuovo regularity of sheaves on the projective space, Enrique Arrondo and Francesco Malaspina obtain an
analogue of Evans and Griffith's theorem (see \cite{Evans-Griff}) on Grassmannians of lines. They give in \cite{Malaspina} two criteria (Theorem 3.1 and 3.2) stating that a vector bundle contains, as a direct summand, an exterior power of the universal sub-bundle or a symmetric power of the universal quotient bundle, if certain cohomologies vanish. They also characterize (Theorem 3.3) those vector
bundles that are direct sums of twists of the above exterior and symmetric powers. 
As a consequence of these results, Alicia Tocino's thesis was carried out under Enrique Arrondo's supervision in 2015, entitled \textit{Cohomological characterization of universal bundles of the Grassmannian of lines}. More concretely, in \cite[Theorem 4.14]{Tocino} a characterization is provided according to directs sums of twists of symmetric powers of the universal quotient bundle over $\G(1,n)$ (also using derived categories and Beilinson's spectral sequence).  

We continue with more vector bundles without intermediate cohomology beyond Grassmannians. Enrique, together with Laura Costa, proves in \cite[Theorem 3.4]{Costa2000} that there exist, up to a twist, only three indecomposable rank-2 bundles without intermediate cohomology on Fano 3-folds. For vector bundles of higher rank, they give a table of concrete examples and characterize which are the Chern classes of vector bundles without
intermediate cohomology (verifying some general conditions), see \cite[Table 4.6 and Theorem 4.9]{Costa2000}. 
With Daniele Faenzi, Enrique focuses on rank-2 vector bundles without intermediate cohomology on prime Fano threefold $X$ of index 1 and genus 12, proving that there are only five different classes and providing a description of their moduli spaces, see \cite{Faenzi}. They mainly use elliptic curves in $X$  and the resolution of the diagonal on $X\times X$ to prove such result. 
In addition, in \cite{Madonna} with Carlo G. Madonna, he studies aCM vector bundles $F$ of rank greater or equal than $3$ on hypersurfaces $X_r$ inside $\p^4$ of degree $r\geq 1$ . The main result is Theorem 1.5, focusing on an aCM vector bundle $F$ of rank $3$ and $4$ on a general quartic threefold $X_4\subset \p^4$ which satisfies a concrete condition, 
defined in Definition 1.2.

\section{Projections and Grassmannian variety}\label{sec-projections}
The study of projective varieties of small codimension that are not  linearly normal, that is, isomorphically  projected from  higher  projective  spaces, is a classical problem. Note that any $n$-dimensional variety can be projected isomorphically to $\p^{2n+1}$, but it produces singular points when it is projected to $\p^{2n}$. So, it is expected that the $n$-dimensional varieties of codimension at most $n$ are linearly normal. In the first place, Francesco Severi proves in \cite{Severi} that the only nondegenerate, that is, not contained in a hyperplane, 
smooth complex surface in $\p^5$ that can be isomorphically projected to $\p^4$ is the Veronese surface. Subsequently, Fyodor L. Zak proves in \cite{Zak} that for $n\geq 2$, the only nondegenerate $n$-dimensional smooth subvariety of $\p^{n(n+3)/2}$ that can be isomorphically projected to $\p^{2n}$ is the double Veronese embedding of $\p^n$, together with a large amount of projectability results in terms of the secant varieties.
Enrique shows in \cite{Arrondo1999Projections} an analogous result for the Grassmannian of lines. He gives a
characterization of the double Veronese embedding of $\p^n$ as the only variety that, under certain general conditions, can be isomorphically projected from $\G(1,2n+1)$ to $\G(1,n+1)$ (Theorem 3.1 and Theorem 3.2).
This was the first step for a huge variety of papers. Also by his own, he studies the same problem for $\G(n-1,n)$ and $\G(n-2,n)$, assuming some general conditions, in \cite[Proposition 2.1 and Theorem 2.3]{Arrondo1998}. 
The topic of projections of Grassmannians is addressed with the thesis, under Enriques's  supervision, of Luca Ugaglia, in 2001, with title \textit{Projection of subvarieties of Grassmannians of lines}. On a similar topic, Jose Carlos Sierra defends  his thesis in 2004, entitled \textit{Proyecciones en Grassmannianas e inmersiones dobles de Veronese}, supervised also by Enrique.
Consequently, in \cite{SierraUgaglia2005} all three of them give a classification of the varieties that are projectable to  $\G(1,n+1)$ coming from  $\G(1,2n)$ and obtaining as a result:

\begin{theo}[\cite{SierraUgaglia2005}, Theorem 3.1]
The only smooth, n-dimensional ($n\geq 3$), nondegenerate,
uncompressed varieties that can be isomorphically projected from $\G(1, 2n)$ to $\G(1, n + 1)$ are the Veronese variety and the blow-up of $\p^n$ in one point.
\end{theo}
\noindent
where $X\subset \G(1,N)$ (${\rm dim}(X)=n$) is said to be \textit{uncompressed} if ${\rm dim}(\Bar{X})=n+1$ with $\Bar{X}$ denoting the union inside $\p^N$ of the lines in $X$ (as defined in Definition 1.3 of \cite{SierraUgaglia2005}).

\noindent
As an extension of these arguments, Enrique, together with Raffaella Paoletti, obtains in \cite{Paoletti} a result concerning Grassmann varieties of higher-dimensional linear subspaces, more concretely, projections from $\G(d-1,nd+d-1)$ into $\G(d-1,n+2d-3)$ under suitable conditions:
\begin{theo}[\cite{Paoletti}, Theorem 3.1]
    Let $X\subset \G(d-1,nd+d-1)$ be a smooth irreducible nondegenerate $n$-dimensional variety such that any two (possibly infinitely close) $(d-1)$-planes of $X$ do not meet. If $X$ has positive defect and is $1$-projectable to $\G(d-1,n+2d-3)$, then $X$ is the Veronese variety.
\end{theo}
\noindent
They also study the relation of this problem with the Steiner bundles over $\p^n$. Indeed, in Proposition 4.3, they prove that the Schwarzenberger bundles never appear if $n\geq 3$. This last remark connects us to the following section.

\section{Steiner and Schwarzenberger bundles}\label{sec-Steiner}
The difficulty of working with a given family of vector bundles is often related with how complicated its resolution is. 
In this direction, a vector bundle $F$ on a projective variety is called a \textit{Steiner bundle} if it is defined as a cokernel of copies of two bundles that form an strongly exceptional pair (see \cite{Miro-Soares}). In particular, in the projective space $\p^n$, the classical definition is given considering a linear resolution of length 1. When the linear map defining a vector bundle is given by a particular cohomological multiplication, the Steiner bundle is called a \textit{Schwarzenberger bundle}. In \cite{Arrondo2010Schwarzenberger}, Enrique introduces a certain class of Steiner bundles that generalize the construction of Schwarzenberger and are therefore called \textit{generalized Schwarzenberger bundles}. Furthermore, he inquires whenever it is possible to describe a given Steiner bundle as a generalized Schwarzenberger bundle. To do this, he defines the concept of jumping subspaces of a Steiner bundle, bounding the dimension of the jumping locus (Theorem 2.8). Finally, he proves that any Steiner bundle whose jumping locus has maximal dimension is in fact a generalized Schwarzenberger bundle and, specifically, every one of such bundles falls in one of four particular families (Theorem 3.7). 

Together with Simone Marchesi, as part of Simone's doctoral project, the latter work has been extended to Grassmannian varieties. These results can be found in Simone's Ph.D. thesis \textit{Jumping pairs of Steiner bundles}, followed by the published work \cite{Marchesi}. More in detail, they gave a general definition of a Steiner bundle on a Grassmannian, finding lower bounds for its possible ranks, and also provide the notion of generalized Schwarzenberger bundle on a Grassmannian. Furthermore, they introduce the notion of jumping pairs associated to a Steiner bundle, bound the dimension of the jumping variety (Theorem 4.9) and prove once again that any Steiner bundle on $\G(k,n)$ whose jumping locus is maximal belongs to a finite list of possible Schwarzenberger bundles (Theorem 5.5). Later on, Enrique, in a joint work with Simone Marchesi and Helena Soares, generalizes the same ideas for smooth projective varieties, see \cite{MarchesiSoares}.

\section{On the path to Hartshorne’s conjecture}\label{sec-Hartshorne}
In 1974, Robin Hartshorne states in \cite{Hart} his celebrated conjecture:
\begin{conj}
    If $X$ is a nonsingular subvariety of dimension $n$ of $\mathbb{P}^N$, and if $n>\frac{2}{3}N$, then $X$ is a complete intersection.
\end{conj}
\noindent
In other words, subvarieties of low codimension must be complete intersections. 
This problem has been deeply investigated by many members of the mathematical community and opened many lines of research, for example the study of curves in a threefold which lead to the definition of reflexive sheaves. Wolf Paul Barth states in \cite{Barth70} that $H^i(X,\mathbb{Q})\cong H^i(\p^N,\mathbb{Q})$ if $i\leq 2n-N$. Subsequently, together with Mogens Esrom Larsen, they prove in \cite{BarthLarse72} that $\pi_1(X)=0$, or equivalently, $X$ is simply connected, if $N\leq 2n-1$. Finally, Mogens Esrom Larsen shows in \cite{Larsen73} that $H^i(X,\Z)\cong H^i(\p^N,\Z)$ if $i\leq 2n-N$. 
In particular, this implies that the Picard group of such $X$ is generated by the class of its hyperplane section if $N \leq  2n-2$ and the cohomology of $X$ ``behaves like the one of a complete intersection''. In the codimension two case, the conjecture can be restated as:
\begin{conj}
    For $N> 6$, any codimension two smooth $X\subset \p^N$ has to be a complete intersection.
\end{conj}
\noindent 
More concretely, as a consequence of the Barth-Larsen theorem, any codimension two variety in the conditions above is subcanonical, that is, the canonical bundle is a multiple of the hyperplane section. In \cite{Arrondo2005Subcanonicity}, Enrique proves the following.
\begin{prop}[\cite{Arrondo2005Subcanonicity}, Proposition 1.1]
Let $X\subset \p^N$ be a smooth subvariety of codimension two. If $N\geq 6$, then $X$ is rationally numerically subcanonical.
\end{prop}
\noindent
He then studies the same problem changing $\p^n$ with a Grassmannian or a quadric. Being the limit and interesting case when the ambient space has dimension six, he focuses on $G(1, 4)$ and $Q_6$, which denotes the smooth six-dimensional quadric, obtaining as main results Theorem 2.1 and 2.2 and Corollary 2.3 for Grassmannians and Theorem 3.2 for quadrics. 

\begin{theo}[\cite{Arrondo2005Subcanonicity}, Theorem 2.1 and Theorem 2.2]
    Any smooth subvariety $X\subset\G(1,n)$ with $n\in\{4,5\}$ of codimension two is rationally numerically subcanonical.
\end{theo}
\begin{coro}[\cite{Arrondo2005Subcanonicity}, Corollary 2.3]
    Let $X\subset \G(1,n)$ be a smooth codimension two subvariety. If $n\geq 4$ then $X$ is rationally numerically subcanonical.
\end{coro}
\begin{theo}[\cite{Arrondo2005Subcanonicity}, Theorem 3.2]
    Let $X\subset Q_6$ be a smooth codimension two subvariety of $Q_6$. Then $X$ is rationally numerically subcanonical if and only if $g_1=g_2$ (where 
    $g_1, g_2$ are the genera of the curves obtained by intersecting $X$ with a three-dimensional linear space of each of the two families of such linear spaces contained in $Q_6$).
\end{theo}

In a joint work with Maria Lucia Fania, see \cite{Fania2006}, he proves that any smooth codimension two projective subvariety of $\G(1,4)$, of degree less than or equal to $25$, is subcanonical, providing a classification of such subvarieties (Theorem 4.1). As a consequence, any smooth codimension two projective subvariety of $\G(1,4)$, which is not of general type, has degree less than or equal to $32$ (Theorem 5.4) and such subvarieties are completely determined in Table 3 of
Example 5.1.

Following this path, Enrique continues to study low-codimension subvarieties in collaboration with Jorge Caravantes, a work which lead to Jorge's Ph.D. thesis: \textit{Sobre el grupo de Picard en subvariedades de codimensión pequeña}. In the corresponding published paper, see \cite{Caravantes}, they try to determine whenever an $n$-dimensional smooth subvariety of an ambient space of dimension at most $2n-2$ inherits the Picard group from the ambient space, and observe that a key step to do this is knowing if the subvariety is simply connected. In particular, they focus on the case of Grassmannian of lines (Theorem 2.1), where some Schubert varieties come into play, and on a product of two projective spaces of the same dimension (Theorem 3.1).

\section{More algebraic geometry topics}\label{sec-various}
This section is a collection, in chronological order, of Enrique's remaining papers dealing with different topics of algebraic geometry that are not easily classifiable due to their variety and complexity.
Together with Manuel Pedreira and his supervisor, Ignacio Sols, he publishes \cite{SolsPedreira1989}, concerning ruled surfaces. They denote by $Q_n\in \p^{n+1}$ the smooth $n$-dimensional quadric and by $H_{d,q}(Q_4)$  the Hilbert scheme of smooth curves of degree $d$ and genus $q$. A ruled surface in $\p^3$ means the image of a ruled surface as a scroll of $\p^3$ with no multiple generators, or equivalently, a smooth curve $C$ in $\G(1,3)$. Moreover, $R_{d,q}(Q_4)$ and $S_{d,q}(Q_4)$ represent the open subschemes of $H_{d,q}(Q_4)$ corresponding respectively to regular and stable ruled surfaces (that is, surfaces in $\p^3$ not having unisecants of degree less than or equal to $d/2$).
Their main result is the following.
\begin{theo}[\cite{SolsPedreira1989}, Theorem]
    If $d\leq 2q+2$ then $R_{d,q}(Q_4)$ and $S_{d,q}(Q_4)$ are irreducible open subschemes of dimension $4d-q+1$ in the same component of $H_{d,q}(Q_4)$.
\end{theo}

Hermann Schubert proves in \cite{Schubert} two formulas  concerning the number of double contacts among the curves of two families in $\p^2$ and also conjectures four other formulas. 
The aim of \cite{Mallavibarrena1990}, written together with Raquel Mallavibarrena and Ignacio Sols, is to give a proof of these six formulas, in the framework of Hilbert's 15th problem, by finding bases of the Chow groups of $\text{Hilb}^2\mathbb{F}$, the Hilbert scheme of the point-line flag variety $\mathbb{F}=\{(P,l)\in \p^2\times \p^{2^*}\,|\, P\in l\}$.

Again with Ignacio Sols, Enrique gives in \cite{Sols1992Bounding}  some bounds on the global sections of vector bundles over a smooth, complete and connected curve and discuss their sharpness. In order to state the main results, we reproduce here the notation they use. Consider $C$ a smooth irreducible curve of genus $g$ and $E$ a rank two vector bundle of degree $d$ on it. Denote by $-e$ the minimum degree of a twist $E\otimes L^{-1}$ having sections, for any line bundle $L$ on $C$ (which is an invariant of the ruled surface $\p(E)$). They propose the following conjecture:
\begin{conj}
    In the above notations, if $-e\leq d\leq 4g-4+e$ and $\p(E)$ is not $C\times \p^1$, then $h^0(E)\leq \frac{d+e}{2}+1$.
\end{conj}
\noindent
They prove that the conjecture is true for the semistable case and for hyperelliptic curves.
\begin{prop}[\cite{Sols1992Bounding}, Proposition 2]
If $C$ is not hyperelliptic, $E$ is semistable and $-e\leq d\leq 4g-4+e$, then $h^0(E)\leq \frac{d}{2}+1$ unless $E$ is either $\mathcal{O}_C\oplus \mathcal{O}_C$ or $\omega_C\oplus \omega_C$.
\end{prop}
\begin{prop}[\cite{Sols1992Bounding}, Proposition 3]
    If $C$ is hyperelliptic and $-e\leq d\leq 4g-4+e$, then $h^0(E)-h^0(\mathcal{E})\leq \frac{d+e}{2}$ and for any values of $d$ and $e$ such that $d\equiv e-2(\text{mod }4)$ there exists a vector bundle $E$ achieving the bound.
\end{prop}
\begin{coro}[\cite{Sols1992Bounding}, Corollary 6]
    Let $E$ be a semistable vector bundle of degree $d$ and rank $R$ which is generically generated by global sections and assume that $h^0(E^\lor \otimes \omega_C)\neq 0$. Then $h^0(E)\leq \frac{d}{2}+R$.
\end{coro}
\noindent
Finally, consider $\mathcal{W}_{d,R}$ the moduli space of semistable rank $R$ vector bundles on $C$ of degree $d$ and $\mathcal{W}_{d,R}^r(C)$ the subscheme parameterizing those bundles $E$ having at least $r+1$ global sections. With these notations, they also prove the following.
\begin{prop}[\cite{Sols1992Bounding}, Proposition 9]
    Let $Y$ be a component of $\mathcal{W}_{d,R}^r(C)$ such that the bundle $E$ corresponding to its generic point is spanned by its sections, has not automorphisms different from the identity and $E\otimes \wedge^R E$ is strongly special. Then, $\dim(Y)\leq (R+1)(\frac{d}{2}-r)$.
\end{prop}

In 1996, his first Ph.D. student, Raquel Díaz, defends her thesis entitled \textit{Matrices de Gram y espacios de ángulos diédricos de poliedros}.

One year later, Enrique with Ignacio Sols and Robert Speiser study in \cite{Sols1997Global} what happens when two embedded varieties, smooth or not (regardless of their dimensions),  make specific contact with each other. The main result is Theorem 7.4 in terms of some data schemes $D_k^r X$ constructed in Section 4.

Aside, jointly with Juana Sendra and  Juan Rafael Sendra, he publishes \cite{Sendras1995Parametric}, in which they extend the classical notion of offset to the concept of generalized offset to a hypersurface, and \cite{Sendras1999} in which the same authors compute the genus of irreducible generalized offset curves to projective irreducible plane curves with only affine ordinary singularities over an algebraically closed field.

Together with Alessandra Bernardi, Enrique publishes \cite{Bernardi2011}. The purpose of this paper is to relate the variety of splitting forms (namely $\text{Split}_d(\p^n)$, see Definition 1.1), that is, the variety whose points are classes of degree $d$ forms splitting as a product of $d$ linear forms in $n+1$ variables with $\G(n-1,n+d-1)$  obtaining also results concerning the higher secant varieties of the varieties of splitting forms. For example, in Theorem 5.4, they provide a result on the intersection between $\G(n-1,n+d-1)$ and $\text{Split}_d(\p^n)$ when $d=3$. The case $d=2$ is also studied in Section 2.

Jointly with Antonio Lanteri and Carla Novelli, he publishes \cite{LanteriNovelli}. They define the notion of ``delta-genus'' for ample vector bundles $E$ of rank two on a smooth projective threefold $X$ as a couple of integers $(\delta_1,\delta_2)$ (Definition 1.1) which extends the classical definition for ample line bundles. 
Furthermore, a classification of $(X, E)$ with low $\delta_1$ and
$\delta_2$ is provided under suitable additional assumptions on $E$. They summarize the main results in the theorem stated in page 138 in the Introduction section.

More recently, with Alessandra Bernardi, Pedro Macias Marques and Bernard Murrain, problems related to skew-symmetric tensor decomposition are considered in \cite{BernardiSkew2021}, but from an algebraic geometry point of view, resulting in the study of higher secant varieties of Grassmannians.
Moreover, from the skew-symmetric action, they define the skew-catalecticant matrices stating the skew-apolarity lemma (Lemma 12) which is the analogue of the classical apolarity lemma for symmetric tensors.

\section{How Enrique understands algebraic geometry}\label{sec-EnriqueAG}
It is common knowledge, especially among his former students, that Enrique has a unique way of understanding mathematics.
Throughout his career, this lead to the publications of papers and extremely useful notes about several topics in algebraic geometry.

For example, we would like to recall \cite{Arrondo2004SerreCorrespondence} and \cite{Arrondo2007Homemade}, where it is possible to find an alternative proof of the Hartshorne-Serre correspondence and represent now a widely used reference for this result. In \cite[Section 1]{Arrondo2007Homemade}, it is recalled the standard approach to the Hartshorne-Serre construction.
\begin{theo}[\cite{Arrondo2007Homemade}, Theorem 1.1: Hartshorne-Serre correspondence]
Let $X$ be a smooth algebraic variety and let $Y$ be a local complete intersection subscheme of codimension two in $X$. Let $N$ be the normal bundle of $Y$ in $X$ and let $L$ be a line bundle on $X$ such that $H^2(X,L^*)=0$. Assume that $\wedge^2N\otimes L^*_{\vert Y}$ has $r-1$ generating global sections $s_1,\ldots, s_{r-1}$. Then there exists a rank $r$ vector bundle $E$ over $X$ such that:
\begin{enumerate}
    \item[(i)] $\wedge^rE=L$;
    \item[(ii)] $E$ has $r-1$ global sections $\alpha_1,\ldots,\alpha_{r-1}$ whose dependency locus is $Y$ and ${s_1 \alpha_1}_{\vert Y}+\ldots+{s_{r-1} \alpha_{r-1}}_{\vert Y}=0$.
\end{enumerate}  Moreover, if $H^1(X,L^*)=0$, conditions $(i)$ and $(ii)$ determine $E$ up to isomorphism.
\end{theo}
\noindent
The elementary proof of this theorem is divided in three sections. In Section 3, he studies the main properties that are required to an open covering of our general ambient variety leading to Lemmas 3.1 and 3.2. 
In order to reproduce here these lemmas, we recall the notation used. Let us take a covering of $Y$ by affine sets $Y\cap U_i$ with $i\in I$ such that $U_i$ is an affine set of $X$, the vector bundle $L$ trivializes on $U_i$ with transition functions $h_{ij}$ and $\mathfrak{J}(U_i)$ is generated by the vanishing of two regular functions $f_i,g_i$ on $U_i$ (where $\mathfrak{J}(U_i)$ denotes the ideal of $Y\cap U_i$ inside $U_i$). Regarding the intersection of two of those open sets, $U_i$, $U_j$, one can find a matrix $A_{ij}$ satisfying 
$\tiny\begin{pmatrix}
    f_i\\ g_i
\end{pmatrix}=A_{ij}\begin{pmatrix}
    f_j\\ g_j
\end{pmatrix}=\begin{pmatrix}
    a_{ij} & b_{ij}\\ c_{ij} & d_{ij}
\end{pmatrix}\begin{pmatrix}
    f_j\\ g_j
\end{pmatrix}$
where $a_{ij}, b_{ij}, c_{ij}, d_{ij}$ are regular functions on $U_i\cap U_j$ and $\det{A_{ij}}$ does not have zeros on $U_i\cap U_j$.
Note that the vector bundle $N$ trivializes on $Y\cap U_i$ and has as transition matrices the restriction $\Bar{A}_{ij}$ of $A_{ij}$ to $Y\cap U_i\cap U_j$. As before, we consider $s_1,\ldots , s_{r-1}$ the global sections generating $\wedge^2N\otimes L^\star$. These sections can be represented locally at each $Y\cap U_i$ by a regular function $\Bar{s}_{it}$ such that $\Bar{s}_{it}=\frac{\det{\Bar{A}_{ij}}}{\Bar{h}_{ij}}\Bar{s}_{jt}$. Since $\Bar{s}_{i1},\ldots, \Bar{s}_{i,r-1}$ do not vanish simultaneously on $Y\cap U_i$, one can refine the covering and assume that there is $t_i\in\{1,\ldots,r-1\}$ such that $\Bar{s}_{i t_i}$ does not have zeros in $Y\cap U_i$. By replacing $U_i$ with its intersection with $\{s_{i t_i}\neq 0\}$, we can assume that $s_{i t_i}$ does not have zeros in $U_i$, so, it is a unit in $\mathcal{O}_X(U_i)$. 
\begin{lem}[\cite{Arrondo2007Homemade}, Lemma 3.1] With the above notations, it is possible to choose regular functions $f_i$, $g_i$ such that $s_{i t_i}=(-1)^{t_i}$. In particular, $\det{\Bar{A}_{ij}}=(-1)^{t_i}\frac{\Bar{h}_{ij}}{\Bar{s}_{j t_i}}.$
\end{lem}
\noindent
The affine covering can be extended to a covering of the whole $X$. So, he covers $X\setminus Y$ by a new affine open sets $U_i$ and defines the matrices $A_{ij}$ for any choice of open sets $U_i, U_j$. More accurately:
\begin{itemize}
    \item If $Y\cap U_i\neq \emptyset \neq Y\cap U_j$,  $A_{ij}=\tiny\begin{pmatrix}
    a_{ij} & b_{ij}\\ c_{ij} & d_{ij}
\end{pmatrix}$, as before.
    \item If $Y\cap U_i= \emptyset = Y\cap U_j$,  $A_{ij}$ is the identity matrix.
    \item If $Y\cap U_i\neq \emptyset = Y\cap U_j$,  $A_{ij}=\tiny\begin{pmatrix}
    u_j & v_j\\ -g_j & f_j
\end{pmatrix}$, with $u_j, v_j$ such that $u_j f_j+v_j f_j=1$.
\item If $Y\cap U_i = \emptyset \neq Y\cap U_j$,  $A_{ij}=\tiny\begin{pmatrix}
    f_i & -v_i\\ g_i & u_i
\end{pmatrix}$, with $u_i, v_i$ such that $u_i f_i+v_i f_i=1$.
\end{itemize}
\begin{lem}[\cite{Arrondo2007Homemade}, Lemma 3.2] With the above choices and notations, it is possible to choose the matrices $A_{ij}$ such that $\det{A_{ij}}=(-1)^{t_i}\frac{h_{ij}}{s_{j t_i}}.$
\end{lem}

\noindent
In Section 4,  he constructs the $r-1$ sections of the desired vector bundle reaching to Lemmas 4.1 and 4.3 and Corollary 4.4. 
For the purpose of stating the first lemma and the corollary, we continue with the appropriate notation. Consider $\alpha_1,\ldots,\alpha_{r-1}$ $r-1$ sections of the vector bundle $E$. Since $\alpha_1,\ldots,\hat{\alpha}_{t_i},\ldots,\alpha_{r-1}$ are linearly independent on $U_i$, it can be extended to a basis of $E_{\vert U_i}$ so that it is possible to represent $\alpha_1,\ldots,\alpha_{r-1}$ on $U_i$ in terms of this basis as the columns of an $r\times (r-2)$ matrix $M_i=\tiny\begin{pmatrix}
    \Delta_{t_i}T'_i\\ T''_i
\end{pmatrix}$, where 
$$T'_i=\tiny\begin{pmatrix}
    1 & 0 & \ldots & -(-1)^{t_i}s_{i1}& \ldots & 0\\
    0 & 1 & \ldots & -(-1)^{t_i}s_{i2} & \ldots & 0\\
    \vdots & & \ddots & \vdots & & \vdots\\
    0 & 0 & \ldots & 1 & \ldots & 0\\
    \vdots & &  & \vdots &\ddots & \vdots\\
    0 & 0 & \ldots & -(-1)^{t_i}s_{i,r-1}& \ldots & 1
\end{pmatrix},  T''_i=\tiny\begin{pmatrix}
    0 & 0 & \ldots & f_i & \ldots & 0\\
    0 & 0 & \ldots & g_i & \ldots & 0
\end{pmatrix}$$
\noindent
and $\Delta_t N$ is the submatrix of $N$ obtained by removing its $t$-th row. Similarly, $N \Delta'_t$ is the submatrix of $N$ obtained by removing its $t$-th column.
\begin{lem}[\cite{Arrondo2007Homemade}, Lemma 4.1]
    For a covering and choices as in Lemma 3.2, if for each $i\in I$ we take $M_i$ as before, then an $r\times r$ matrix $Z_{ij}=\tiny\begin{pmatrix}
        P_{ij} & Q_{ij}\\ R_{ij} & S_{ij}
    \end{pmatrix}$ satisfies the equality $M_i=Z_{ij}M_j$ if and only if the following equalities hold:
    \begin{enumerate}
        \item $P_{ij}=\Delta_{t_i}T'_i\Delta'_{t_j}$
        \item $R_{ij}=T''_i\Delta'_{t_j}$
        \item $Q_{ij}{\tiny \begin{pmatrix}
            f_j\\g_j
        \end{pmatrix}}=(-1)^{t_j}\Delta_{t_i}T'_i{\tiny \begin{pmatrix}
            s_{j1}\\\vdots\\s_{j r-1}
        \end{pmatrix}}$
        \item $S_{ij}{\tiny \begin{pmatrix}
            f_j\\g_j
        \end{pmatrix}}=(-1)^{t_j}s_{j t_i}{\tiny \begin{pmatrix}
            f_i\\g_i
        \end{pmatrix}}$
    \end{enumerate}
    Moreover, such a matrix always exists and, when taking $A_{ij}$ as in Lemma 3.2, it follows $\det{S_{ij}}=(-1)^{t_i}s_{jt_i}h_{ij}$ and $\det{Z_{ij}}=h_{ij}$.
\end{lem}
\begin{coro}[\cite{Arrondo2007Homemade}, Corollary 4.4]
    If the matrices $\{Z_{ij}\}_{i,j\in I}$ are chosen as in the previous lemma, then for any $i,j,k\in I$ there exist regular functiones $\beta_{ijk1},\ldots, \beta_{ijk,r-1}$ on $U_i\cap U_j\cap U_k$ such that $Z_{ik}-Z_{ij}Z_{jk}=\begin{pmatrix}
        0 & B_{ijk}
    \end{pmatrix}$, with
    $$B_{ijk}=\begin{pmatrix}
        Q_{ik}-P_{ij}Q_{jk}-Q_{ij}S_{jk}\\
        S_{ik}-R_{ij}Q_{jk}-S_{ij}S_{jk}
    \end{pmatrix}=\begin{pmatrix}
        \beta_{ijk1}\\\vdots\\ \hat{\beta}_{ijkt_i} \\ \vdots \\ \beta_{ijk,r-1} \\ \beta_{ijkt_i}f_i \\ \beta_{ijkt_i}g_i
    \end{pmatrix} \begin{pmatrix}
        g_k, -f_k
    \end{pmatrix}.$$
\end{coro}

\noindent
The end of the proof culminates with Proposition 5.6 in Section 5 in which the existence and uniqueness of the vector bundle is proven due to cohomological conditions on $L^*$.
So, the correspondence is obtained gluing together the explicit local description in order to construct the section of a vector bundle. 

Moreover, he provides \textit{Another elementary proof of the Nullstellensatz} in \cite{Arrondo2006Nullstellensatz} and \textit{The Nullstellensatz without the Axiom of Choice} in \cite{arrondo2021nulls}. In addition, he gives a new innovative approach to construct representations of finite groups without the necessity of calculating character tables in \cite{arrondo2021representation} under the name \textit{Representation theory of finite groups through (basic) algebraic geometry}. 

Finally, we would like to remind the reader of Enrique's useful lecture notes on several topics, starting from lecture notes on bachelor degree level to arrive to more complicated aspects of algebraic geometry. In particular, we would like recall \cite{arrondo1996notesGrassmannians}, which provides an introduction to Grassmannians and its subvarieties, that, as we could aprreciate in this survey, represent a common factor in Enrique's work. All of these notes are available at Enrique's webpage \url{https://blogs.mat.ucm.es/arrondo}. 

\bigskip
\textbf{Acknowledgements:} We thank the anonymous referee for his/her valuable suggestions that improved the quality of this paper.

\bibliographystyle{plain}
\bibliography{ref}

\begin{thebibliography}{10}

\bibitem{arrondo1996notesGrassmannians}
Enrique Arrondo.
\newblock {\em Subvarieties of Grassmannians}, volume~10.
\newblock Lecture Note Series Dipartimento di Matematica Univ. Trento, 1996.

\bibitem{Arrondo1998}
Enrique Arrondo.
\newblock The universal {$\text{rank-}(n-1)$} bundle on {$G(1,n)$} restricted to subvarieties.
\newblock {\em Collect. Math.}, 49(2-3):173--183, 1998.
\newblock Dedicated to the memory of Fernando Serrano.

\bibitem{Arrondo1999Projections}
Enrique Arrondo.
\newblock Projections of {G}rassmannians of lines and characterization of {V}eronese varieties.
\newblock {\em J. Algebraic Geom.}, 8(1):85--101, 1999.

\bibitem{Arrondo2002LineCongr}
Enrique Arrondo.
\newblock Line congruences of low order.
\newblock {\em Milan J. Math.}, 70:223--243, 2002.

\bibitem{Arrondo2004SerreCorrespondence}
Enrique Arrondo.
\newblock Deriving the {S}erre correspondence by hand.
\newblock In {\em Mathematical contributions in honor of {P}rofessor {E}nrique {O}uterelo {D}om\'{\i}nguez ({S}panish)}, Homen. Univ. Complut., pages 61--72. Editorial Complutense, Madrid, 2004.

\bibitem{Arrondo2005Subcanonicity}
Enrique Arrondo.
\newblock Subcanonicity of codimension two subvarieties.
\newblock {\em Rev. Mat. Complut.}, 18(1):69--80, 2005.

\bibitem{Arrondo2006Nullstellensatz}
Enrique Arrondo.
\newblock Another elementary proof of the {N}ullstellensatz.
\newblock {\em Amer. Math. Monthly}, 113(2):169--171, 2006.

\bibitem{Arrondo2007Homemade}
Enrique Arrondo.
\newblock A home-made {H}artshorne-{S}erre correspondence.
\newblock {\em Rev. Mat. Complut.}, 20(2):423--443, 2007.

\bibitem{Arrondo2010Schwarzenberger}
Enrique Arrondo.
\newblock Schwarzenberger bundles of arbitrary rank on the projective space.
\newblock {\em J. Lond. Math. Soc. (2)}, 82(3):697--716, 2010.

\bibitem{arrondo2021nulls}
Enrique Arrondo.
\newblock The nullstellensatz without the axiom of choice.
\newblock {\em arXiv preprint}, 2020.

\bibitem{arrondo2021representation}
Enrique Arrondo.
\newblock Representation theory of finite groups through (basic) algebraic geometry.
\newblock {\em arXiv preprint}, 2021.

\bibitem{Bernardi2011}
Enrique Arrondo and Alessandra Bernardi.
\newblock On the variety parameterizing completely decomposable polynomials.
\newblock {\em J. Pure Appl. Algebra}, 215(3):201--220, 2011.

\bibitem{BernardiSkew2021}
Enrique Arrondo, Alessandra Bernardi, Pedro~Macias Marques, and Bernard Mourrain.
\newblock Skew-symmetric tensor decomposition.
\newblock {\em Commun. Contemp. Math.}, 23(2):Paper No. 1950061, 29, 2021.

\bibitem{BertoliniTurrini1994Classif}
Enrique Arrondo, Marina Bertolini, and Cristina Turrini.
\newblock Classification of smooth congruences with a fundamental curve.
\newblock In {\em Projective geometry with applications}, volume 166 of {\em Lecture Notes in Pure and Appl. Math.}, pages 43--56. Dekker, New York, 1994.

\bibitem{BertoliniTurrini1998SmallDegree}
Enrique Arrondo, Marina Bertolini, and Cristina Turrini.
\newblock Congruences of small degree in {$G(1,4)$}.
\newblock {\em Comm. Algebra}, 26(10):3249--3266, 1998.

\bibitem{BertoliniTurrini2000Grass}
Enrique Arrondo, Marina Bertolini, and Cristina Turrini.
\newblock Quadric bundle congruences in {$G(1,n)$}.
\newblock {\em Forum Math.}, 12(6):649--666, 2000.

\bibitem{BertoliniTurrini2001Surfaces}
Enrique Arrondo, Marina Bertolini, and Cristina Turrini.
\newblock A focus on focal surfaces.
\newblock {\em Asian J. Math.}, 5(3):535--560, 2001.

\bibitem{BertoliniTurrini2005Loci}
Enrique Arrondo, Marina Bertolini, and Cristina Turrini.
\newblock Focal loci in {$G(1,N)$}.
\newblock {\em Asian J. Math.}, 9(4):449--472, 2005.

\bibitem{BertoliniTurrini2011}
Enrique Arrondo, Marina Bertolini, and Cristina Turrini.
\newblock On the ampleness of the normal bundle of line congruences.
\newblock {\em Forum Math.}, 23(2):223--244, 2011.

\bibitem{Caravantes}
Enrique Arrondo and Jorge Caravantes.
\newblock On the {P}icard group of low-codimension subvarieties.
\newblock {\em Indiana Univ. Math. J.}, 58(3):1023--1050, 2009.

\bibitem{Cobo}
Enrique Arrondo and Sof\'{\i}a Cobo.
\newblock On the stability of the universal quotient bundle restricted to congruences of low degree of {$\Bbb G(1,3)$}.
\newblock {\em Ann. Sc. Norm. Super. Pisa Cl. Sci. (5)}, 9(3):503--522, 2010.

\bibitem{Costa2000}
Enrique Arrondo and Laura Costa.
\newblock Vector bundles on {F}ano 3-folds without intermediate cohomology.
\newblock {\em Comm. Algebra}, 28(8):3899--3911, 2000.

\bibitem{Faenzi}
Enrique Arrondo and Daniele Faenzi.
\newblock Vector bundles with no intermediate cohomology on {F}ano threefolds of type {$V_{22}$}.
\newblock {\em Pacific J. Math.}, 225(2):201--220, 2006.

\bibitem{Fania2006}
Enrique Arrondo and Maria~Lucia Fania.
\newblock Evidence to subcanonicity of codimension two subvarieties of {$\Bbb G(1,4)$}.
\newblock {\em Internat. J. Math.}, 17(2):157--168, 2006.

\bibitem{BeatrizGrana1999}
Enrique Arrondo and Beatriz Gra\~{n}a.
\newblock Vector bundles on {$G(1,4)$} without intermediate cohomology.
\newblock {\em J. Algebra}, 214(1):128--142, 1999.

\bibitem{BeatrizGrana2006}
Enrique Arrondo and Beatriz Gra\~{n}a Otero.
\newblock Congruences on {$G(1,4)$} with split universal quotient bundle.
\newblock {\em J. Geom. Phys.}, 56(6):1057--1067, 2006.

\bibitem{Mark1993}
Enrique Arrondo and Mark Gross.
\newblock On smooth surfaces in {${\rm Gr}(1,\bold P^3)$} with a fundamental curve.
\newblock {\em Manuscripta Math.}, 79(3-4):283--298, 1993.

\bibitem{LanteriNovelli}
Enrique Arrondo, Antonio Lanteri, and Carla Novelli.
\newblock A notion of {$\Delta$}-multigenus for certain rank two ample vector bundles.
\newblock {\em Kodai Math. J.}, 36(1):137--153, 2013.

\bibitem{Madonna}
Enrique Arrondo and Carlo~G. Madonna.
\newblock Curves and vector bundles on quartic threefolds.
\newblock {\em J. Korean Math. Soc.}, 46(3):589--607, 2009.

\bibitem{Malaspina}
Enrique Arrondo and Francesco Malaspina.
\newblock Cohomological characterization of vector bundles on {G}rassmannians of lines.
\newblock {\em J. Algebra}, 323(4):1098--1106, 2010.

\bibitem{Mallavibarrena1990}
Enrique Arrondo, Raquel Mallavibarrena, and Ignacio Sols.
\newblock Proof of {S}chubert's conjectures on double contacts.
\newblock In {\em Enumerative geometry ({S}itges, 1987)}, volume 1436 of {\em Lecture Notes in Math.}, pages 1--29. Springer, Berlin, 1990.

\bibitem{Marchesi}
Enrique Arrondo and Simone Marchesi.
\newblock Jumping pairs of {S}teiner bundles.
\newblock {\em Forum Math.}, 27(6):3233--3267, 2015.

\bibitem{MarchesiSoares}
Enrique Arrondo, Simone Marchesi, and Helena Soares.
\newblock Schwarzenberger bundles on smooth projective varieties.
\newblock {\em J. Pure Appl. Algebra}, 220(9):3307--3326, 2016.

\bibitem{Paoletti}
Enrique Arrondo and Raffaella Paoletti.
\newblock Characterization of {V}eronese varieties via projection in {G}rassmannians.
\newblock In {\em Projective varieties with unexpected properties}, pages 1--12. Walter de Gruyter, Berlin, 2005.

\bibitem{SolsPedreira1989}
Enrique Arrondo, Manuel Pedreira, and Ignacio Sols.
\newblock On regular and stable ruled surfaces in {${\bf P}^3$}.
\newblock In {\em Algebraic curves and projective geometry ({T}rento, 1988)}, volume 1389 of {\em Lecture Notes in Math.}, pages 1--15. Springer, Berlin, 1989.

\bibitem{Sendras1995Parametric}
Enrique Arrondo, Juana Sendra, and J.~Rafael Sendra.
\newblock Parametric generalized offsets to hypersurfaces.
\newblock volume~23, pages 267--285. 1997.
\newblock Parametric algebraic curves and applications (Albuquerque, NM, 1995).

\bibitem{Sendras1999}
Enrique Arrondo, Juana Sendra, and J.~Rafael Sendra.
\newblock Genus formula for generalized offset curves.
\newblock {\em J. Pure Appl. Algebra}, 136(3):199--209, 1999.

\bibitem{SierraUgaglia2005}
Enrique Arrondo, Jos\'{e}~Carlos Sierra, and Luca Ugaglia.
\newblock Classification of {$n$}-dimensional subvarieties of {$G(1,2n)$} that can be projected to {$G(1,n+1)$}.
\newblock {\em Bull. London Math. Soc.}, 37(5):673--682, 2005.

\bibitem{Sols1989}
Enrique Arrondo and Ignacio Sols.
\newblock Classification of smooth congruences of low degree.
\newblock {\em J. Reine Angew. Math.}, 393:199--219, 1989.

\bibitem{Sols1992Bounding}
Enrique Arrondo and Ignacio Sols.
\newblock Bounding sections of bundles on curves.
\newblock In {\em Complex projective geometry ({T}rieste, 1989/{B}ergen, 1989)}, volume 179 of {\em London Math. Soc. Lecture Note Ser.}, pages 24--31. Cambridge Univ. Press, Cambridge, 1992.

\bibitem{Sols1992Congruences}
Enrique Arrondo and Ignacio Sols.
\newblock On congruences of lines in the projective space.
\newblock {\em M\'{e}m. Soc. Math. France (N.S.)}, (50):96, 1992.

\bibitem{Sols1997Global}
Enrique Arrondo, Ignacio Sols, and Robert Speiser.
\newblock Global moduli for contacts.
\newblock {\em Ark. Mat.}, 35(1):1--57, 1997.

\bibitem{Tocino}
Enrique Arrondo and Alicia Tocino.
\newblock Cohomological characterization of universal bundles of {$\Bbb{G}(1,n)$}.
\newblock {\em J. Algebra}, 540:206--233, 2019.

\bibitem{Barth70}
Wolf~Paul Barth.
\newblock Transplanting cohomology classes in complex-projective space.
\newblock {\em Amer. J. Math.}, 92:951--967, 1970.

\bibitem{BarthLarse72}
Wolf~Paul Barth and Mogens~Esrom Larsen.
\newblock On the homotopy groups of complex projective algebraic manifolds.
\newblock {\em Math. Scand.}, 30:88--94, 1972.

\bibitem{GBS}
R.-O. Buchweitz, G.-M. Greuel, and F.-O. Schreyer.
\newblock Cohen-{M}acaulay modules on hypersurface singularities. {II}.
\newblock {\em Invent. Math.}, 88(1):165--182, 1987.

\bibitem{Peskine}
Geir Ellingsrud and Christian Peskine.
\newblock Sur les surfaces lisses de {${\bf P}_4$}.
\newblock {\em Invent. Math.}, 95(1):1--11, 1989.

\bibitem{Evans-Griff}
E.~Graham Evans and Phillip Griffith.
\newblock The syzygy problem.
\newblock {\em Ann. of Math. (2)}, 114(2):323--333, 1981.

\bibitem{Livorni}
Maria~Lucia Fania and Elvira~Laura Livorni.
\newblock Degree nine manifolds of dimension greater than or equal to {$3$}.
\newblock {\em Math. Nachr.}, 169:117--134, 1994.

\bibitem{Fano}
Gino Fano.
\newblock Nuove ricerche sulle congruenze di rette del $3^\text{o}$ ordine prive di linea singolare.
\newblock {\em Atti R. Acc. Sci. Torino}, 51:1--79, 1901.

\bibitem{Giraldo}
Luis Giraldo.
\newblock A bound for the arithmetic genus of curves in {G}rassmannians.
\newblock {\em Forum Math.}, 12(6):667--669, 2000.

\bibitem{Hart}
Robin Hartshorne.
\newblock Varieties of small codimension in projective space.
\newblock {\em Bull. Amer. Math. Soc.}, 80:1017--1032, 1974.

\bibitem{Horrocks}
Geoffrey Horrocks.
\newblock Vector bundles on the punctured spectrum of a local ring.
\newblock {\em Proc. London Math. Soc. (3)}, 14:689--713, 1964.

\bibitem{Ionescu}
Paltin Ionescu.
\newblock Embedded projective varieties of small invariants. {III}.
\newblock In {\em Algebraic geometry ({L}'{A}quila, 1988)}, volume 1417 of {\em Lecture Notes in Math.}, pages 138--154. Springer, Berlin, 1990.

\bibitem{Knorrer}
Horst Kn\"{o}rrer.
\newblock Cohen-{M}acaulay modules on hypersurface singularities. {I}.
\newblock {\em Invent. Math.}, 88(1):153--164, 1987.

\bibitem{Larsen73}
Mogens~Esrom Larsen.
\newblock On the topology of complex projective manifolds.
\newblock {\em Invent. Math.}, 19:251--260, 1973.

\bibitem{Miro-Soares}
Rosa~Maria Mir\'{o}-Roig and Helena Soares.
\newblock Cohomological characterisation of {S}teiner bundles.
\newblock {\em Forum Math.}, 21(5):871--891, 2009.

\bibitem{OttavianiGrass2}
Giorgio Ottaviani.
\newblock Crit\`eres de scindage pour les fibr\'{e}s vectoriels sur les {G}rassmanniennes et les quadriques.
\newblock {\em C. R. Acad. Sci. Paris S\'{e}r. I Math.}, 305(6):257--260, 1987.

\bibitem{OttavianiGrass}
Giorgio Ottaviani.
\newblock Some extensions of {H}orrocks criterion to vector bundles on {G}rassmannians and quadrics.
\newblock {\em Ann. Mat. Pura Appl. (4)}, 155:317--341, 1989.

\bibitem{Schubert}
Hermann Schubert.
\newblock Anzahlgeometrische behandlung des dreiecks.
\newblock {\em Math. Ann.}, 17:153--212, 1880.

\bibitem{Severi}
Francesco Severi.
\newblock Intorno ai punti doppi impropri di una superficie generale dello spazio a quattro dimensioni, e a’suoi punti tripli apparenti.
\newblock {\em Rend. Circ. Matem. Palermo}, 15:33--51, 1901.

\bibitem{Verra}
Alessandro Verra.
\newblock Smooth surfaces of degree {$9$} in {$G(1,3)$}.
\newblock {\em Manuscripta Math.}, 62(4):417--435, 1988.

\bibitem{Zak}
Fyodor~L. Zak.
\newblock {\em Tangents and secants of algebraic varieties}, volume 127 of {\em Translations of Mathematical Monographs}.
\newblock American Mathematical Society, Providence, RI, 1993.
\newblock Translated from the Russian manuscript by the author.

\end{thebibliography}
\end{document}